## DISCUSSION

By Linglong Kong and Ivan Mizera[1]

*University of Alberta*

**1. Introduction.** Hallin, Paindaveine and Šiman—hereafter HPŠ—are to be congratulated for the appreciation their paper is receiving from *The Annals of Statistics*. We personally are indebted for the attention they devoted to our paper, [11]—hereafter KM. The topic is quite delicate and multi-faceted, hence all the misunderstandings we allege below show deficiencies of our exposition rather than anything else.

**2. This is not a quantile.** The altered title of Magritte's painting entitling this section can be seen as a gnomic expression of some subtle, nevertheless fundamental differences between HPŠ and KM, a distinct "philosophy" in the language of HPŠ. While HPŠ state in their first sentence their intent to "propose a definition of multivariate quantiles," KM in their second sentence maintain

> While such an objective could be mistaken for yet another attempt in the ongoing quest for "multivariate quantiles"... we would like to stress that we differ in the position that *no multivariate generalization of the quantile concept may be needed at all*—...

Instead, our intention was to discuss "certain aspects of using quantiles to obtain insights about multivariate data," suggesting that interesting insights about multivariate data can be obtained by looking at *univariate* quantiles of projections. Nothing beyond concepts already well established in data-analytic practice: see [12] regarding univariate quantiles (the concept we feel primarily rooted in the order of the numerical scale rather than in their secondary $L^1$ characterization), and any textbook on multivariate analysis regarding projections.

We admit that KM might have not accentuated the word "univariate" enough—nevertheless, nothing like "projection quantiles" is mentioned, nor

Received August 2009.

[1]Supported by the Natural Sciences and Engineering Research Council of Canada.

anything like $\tau$-quantile or $\mathbf{q}_{\mathrm{KM};\nu\mathbf{u}}$ is introduced there. Put verbally, to escape sometimes serpentine notation: KM *do not consider lines orthogonal to the quantile directions "multivariate quantiles,"* but merely consider them a graphical tool to visualize *univariate* quantiles *of* the corresponding *projections.* If HPŠ write

> "The resulting quantile contours (the collections, for fixed $\tau$, of $\mathbf{q}_{\mathrm{KM};\nu\mathbf{u}}$'s do not enjoy the properties (independence with respect to to the choice of an origin, affine-equivariance, nestedness, etc.) one is expecting from a quantile concept."

We would like to repeat that "quantile biplots" (which HPŠ use as a straw man in their Section 3.2 and elsewhere) were brought up by KM *not as a recommendation for practical use*, but as something that shows the pitfalls of attempts to summarize quantiles of projections graphically, something that eventually vindicates the use of depth contours—which is what HPŠ arrive to anyway, if it comes to *quantile contours.* If HPŠ surmise

> Contrary to Kong and Mizera's, however, the $\pi_{\tau\mathbf{u}}$ quantile hyperplanes do enjoy all the desirable properties of a well-behaved quantile concept.

We ask: contrary to Kong and Mizera's what? Contrary to "quantile contours"—which, being *collections* of objects are compared here to "quantile hyperplanes," that is, *single* objects? Or contrary to something else KM threw in the garbage? There is no need to construct an artificial conflict with an imaginary adversary: if HPŠ really want to compare multivariate quantiles, they can take any proposal from the long string of references given in their second paragraph, the string that (rightly!) does not contain KM.

To dispel an impression that we are engaging here in some nihilistic wordplay, let us consider an illustrative example. It will be a comparison; *not of quantile concepts*, but of possible ways of obtaining some interpretable insights from multivariate data. Suppose that we have a dataset concerning weight and height of several subjects. A possible quantity to look at is the Body Mass Index, BMI, as defined in Section 7 of HPŠ. If we use the trick we learned from [15], and consider the data on the logarithmic scale, then $\log(\mathrm{BMI})$ is $\log(\mathrm{weight}) - 2\log(\mathrm{height})$, a linear combination of the new variables; its (univariate!) quantiles correspond to those of the projection into the "log(BMI) direction," the direction of the vector $(1,-2)$. To visualize this, we may rotate the data so that the "log(BMI) direction" becomes vertical (Figure 1). The values of log(BMI) (up to a scale factor, which we do not worry about, because we are after univariate quantiles which are scale-equivariant) then correspond to the projections of the datapoints onto the ordinate; a horizontal line drawn through the chosen quantile ($p = 0.75$, say) then indicates which subjects have (log) BMI above the third quartile, which below, and how the data are divided in this respect—not a big deal perhaps, but something real, something that, say, a lay person like a health professional may be able to understand.

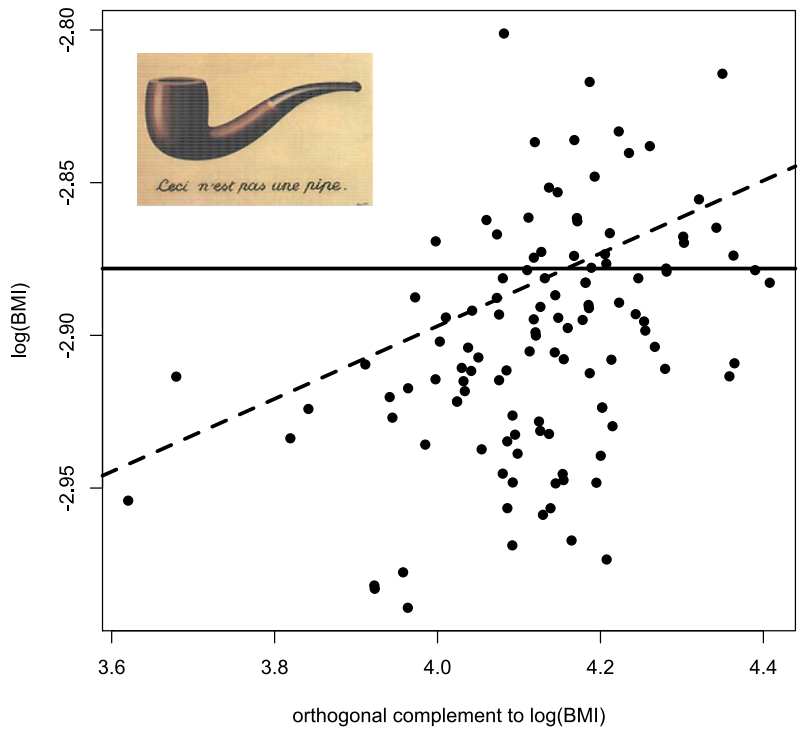

Fig. 1. *The solid line is not a quantile, but merely indicates the location of the third sample quartile of the projection of the data to the ordinate. The broken line is the proposed multivariate quantile. Top left corner is Magritte.*

Now, given the same direction and same $p = 0.75$, HPŠ would suggest to fit the regression quantile instead (the broken line in Figure 1). Let us attempt an interpretation. We know that the line shows the (estimated) *conditional* (0.75)-quantile. This sounds promising—however: conditional on *what*? The projection of the data to the direction orthogonal to the log(BMI), showing a quantity proportional to $2\log(\text{weight}) + \log(\text{height})$? What is the meaning of the latter? Even if we take an affine transformation instead (relying on equivariance), and transform the data so that the horizontal axis corresponds to some interpretable variable, for instance to log(weight): what is then the meaning of the quantiles of log(BMI) *conditionally* on log(weight)? Can we explain that—to a lay person, or to anybody?

**3. It's better to solve the right problem approximately.** ...than the wrong problem exactly. The ubiquitous quote attributed to Tukey by [1] is invoked here not because we would like to suggest that we solved the "right" and somebody else the "wrong" problem, but because we believe that it expresses the principal tenet of practical data analysis: the emphasis is on having the job done, and much less on finding the "best" way of doing it.

Our principal goal in KM—and here we probably erred again not stating that more bluntly—was to show how depth contours can be "regressed on covariates." HPŠ comment that this issue is "briefly addressed" in our Section 11.3; well, other issues might take more space, but we believe that the solution we gave is complete, including a numerical realization and also theoretical analysis (in Section 9 of KM).

After seeing an early presentation of what later became [15], we wondered whether a similar construction cannot be accomplished using depth contours, as earlier suggested in the growth-chart context by [14]. It may be of interest here that our initial idea was to use the Laine's characterization of depth contours, as published in [8]—that is, along the lines how perhaps HPŠ would suggest to tackle the problem. Our later decision to abandon this approach may be interpreted as our ineptitude, or lack of good judgment—although we do not feel it exactly that way. Some technical aspects may be of interest.

Abstracting from certain computational details, the problem is the following: for any given direction $u \in \mathbb{S}^1$, consider the response values $(y_k^u, z_k^u)$, arising from rotating the original bivariate response $(y_k, z_k)$ so that $y_k^u$, the transformed $y_k$, corresponds to the ordinate. What is needed then is a quantile regression of $y_k^u$ on $z_k^u$, and also on the time covariate $x_k$, the regression that satisfies the following: the dependence of $y_k^u$ on $z_k^u$, for fixed $x_k$, should be *linear* (so that Laine's characterization applies), but at the same time the dependence of $y_k^u$ on $x_k$ (with $z_k^u$ fixed, say) should have a flexible, that is, *nonparametric* form. The latter is important, because after seeing several real biometric growth charts—see KM or [15]—it is clear that hardly any parametric model can be expected to describe those data well.

The need to fit a parametric–nonparametric quantile regression makes the problem intricate, and not that straightforward within certain methodologies (think kernel regression, for instance). Fitting a model where everything is linear would be easy (as can be seen from HPŠ); however, the nonparametric dependence on $x_k$ is really vital for achieving a good description of the data here. The only strategy *we* were able to conceive was via a penalized approach adapted to our particular situation. Let us explain its main features; for simplicity, the superscript $u$ will be suppressed in the notation.

The standard formulation for bivariate nonparametric fitting—as discussed, for instance, in [10]—seeks the fit $f$ as a solution of the optimization task

$$(3.1) \qquad \sum_{k=1}^{n} \rho_\tau(y_k - f(z_k, x_k)) + \lambda J(f) \hookrightarrow \min_f!,$$

where $J$ can be, depending on what amount of smoothness we require from the fit, either the classical thin-plate penalty

$$(3.2) \qquad J(f) = \iint \left(\frac{\partial^2 f}{\partial z^2}\right)^2 + \left(\frac{\partial^2 f}{\partial z\,\partial x}\right)^2 + \left(\frac{\partial^2 f}{\partial x^2}\right)^2 dz\,dx$$

or some of its $L^1$ variants, say,

$$(3.3) \qquad J(f) = \iint \left|\frac{\partial^2 f}{\partial z^2}\right| + \left|\frac{\partial^2 f}{\partial z\,\partial x}\right| + \left|\frac{\partial^2 f}{\partial x^2}\right| dz\,dx.$$

It is well-known—see, for instance, [6]—that (3.1) is the Lagrangian version of the constrained regularization formulation

$$(3.4) \qquad \sum_{k=1}^{n} \rho_\tau(y_k - f(z_k, x_k)) \hookrightarrow \min_f! \qquad \text{subject to } J(f) \leq \Lambda,$$

as both yield the same solutions for appropriately related values of $\lambda$ and $\Lambda$. The penalties (3.2) or (3.3), however, yield $f$ that is nonparametric (that is, generally nonlinear) in *both* $z$ and $x$. A possible way to make $f$ *linear* in $z$ for fixed $x$ is to enforce $\partial^2 f/\partial z^2$ to be zero; the formulation (3.4) then becomes

$$\sum_{k=1}^{n} \rho_\tau(y_k - f(z_k, x_k)) \hookrightarrow \min_f! \qquad \text{subject to } \left|\frac{\partial^2 f}{\partial z^2}\right| \leq 0 \text{ and } \tilde{J}(f) \leq \Lambda,$$
(3.5)
where $\tilde{J}$ is a penalty obtained, say, from (3.2) or (3.3) by dropping the $\partial^2 f/\partial z^2$ term from the integral. For instance, (3.3) leads to the penalty

$$(3.6) \qquad \tilde{J}(f) = \iint \left|\frac{\partial^2 f}{\partial z\,\partial x}\right| + \left|\frac{\partial^2 f}{\partial x^2}\right| dz\,dx;$$

however, it is also possible to consider the penalty

$$(3.7) \qquad \tilde{J}(f) = \iint \left|\frac{\partial^2 f}{\partial x^2}\right| dz\,dx$$

instead. It is not clear whether the mixed-derivative term has to be kept or dropped—whether the dependence on $z$ and $x$ is additive or not. The Lagrange equivalent of (3.5) is a formulation with two tuning parameters,

$$(3.8) \qquad \sum_{k=1}^{n} \rho_\tau(y_k - f(z_k, x_k)) + \lambda_1 \left|\frac{\partial^2 f}{\partial z^2}\right| + \lambda_2 \tilde{J}(f) \hookrightarrow \min_f!$$

with $\lambda_1$ set very large, so that it effectively pushes the derivative in $z$ to zero, and $\lambda_2$ chosen to obtain a fit flexible in $x$. For the quadratic penalty (3.2), the alternatives are analogous; the use of the squared second derivative with respect to (3.5) or (3.8) amounts to the same effect, the only substantial difference is the form of $\tilde{J}$.

Although we were able to implement a solution in this vein, and the numerical experiments gave somewhat promising results, there were serious issues. First, it is not clear whether to use (3.6) or (3.7)—both gave plausible, albeit somewhat different fits. Second, there is the notorious problem of selection of $\lambda_2$—not that unsurmountable for *one* regression, but here we have to deal with a parametrized family of regressions, an interesting conceptual problem, which we would, however, prefer to tackle separately at some other time. Third, and most important: we suddenly realized that we are chasing two rabbits simultaneously—even if we succeeded in overcoming the technological issues, any criticism of the methodology we would adopt for quantile regression would negatively affect also the proposed methodology of regressing on depth contours: in other words, if the results would not be convincing, it would not be clear whether it was problem of the directional characterization of depth contours, or the problem of the nonparametric quantile regression methodology.

In such an unfortunate situation the solution came as the proverbial egg of Columbus: why not fit constants instead of lines?! That amounts to projecting $(z_k, y_k)$ into the direction of $u$, which yields a *univariate* response $w_k^u$—which is then fitted by nonparametric quantile regression on $x_k$. We sacrifice some nice mathematics (curiously, HPŠ find the property we use "somewhat surprising"—but it is in fact obvious; it is Laine's characterization that is nontrivial!), but we gain a lot of autonomy: the required nonparametric quantile regressions now involve only one covariate, which means that there are quite a few methods out there proposed for this task. The users can choose one they prefer—for instance, if they want to avoid "quantile crossing phenomenon," presented as unavoidable in Section 6 of HPŠ, and thus obtain the contours that *are* nested, they can use recent methods dealing with this phenomenon, proposed, for instance, by [4] or [5].

So, by switching to univariate projections, we (i) got the job done; (ii) obtained approach that is independent of the finesses of the chosen quantile regression methodology; and (iii) had our attention turned to the data-analytic insights obtained by univariate projections, the philosophy that enabled us to avoid pursuit of nebulous concepts like multivariate quantiles.

If we also took some satisfaction in making one step ahead (or sidewise) of Laine, we would like to stress that this in no case meant the universal a priori dismissal of that (or any other) idea as unsuitable to lead to the solution of the same problem—and perhaps better and even easier. Provided, of course, it *is* the same problem—we maintain that the nonparametric dependence on $x$ is really essential. HPŠ hint at the possible application of local linear methods here; we will await their, or anybody else's relevant implementation with great curiosity. As put by Anonymous (quoted with permission): "The early bird gets the worm. The second mouse gets the cheese."

**4. One country, two systems.** Despite all the subtle differences aired here, we agree that HPŠ and KM have a lot in common: they both address depth contours by adopting a sort of directional approach. Regression aspects having been discussed to some extent in the previous section, we concentrate now on the possible impact on the computation of the depth contours (in the "static" situation, without covariates).

Here HPŠ should be credited for striking upon an important aspect: exact algorithms for higher dimensions—which do exist (despite the impression of HPŠ); see [2, 3] and the references there—indeed involve elements of linear programming. The topic is considerably challenging: the classical result of [7] shows that the problem is NP-hard when also the dimension is the part of the input, and gives lower bounds on best possible complexities of *exact* algorithms. The theoretical complexities derived by [2] depend themselves on the complexity of linear programming—a well-known chapter for itself.

However, rather than speculating on *potential* higher-dimensional implementations, let us comment on our personal experience with good old dimension two. Although obtaining a working routine for depth contours seems to be straightforward there, given the existing literature and (much less abundant than the latter, but still existing) implementations, the truth is not that simple. The implementation of some of the proposed algorithms, although in principle possible, would be so laborious that from the practical point of view, these algorithms may be considered unimplementable. Other algorithms may work reasonably well for moderate sample sizes, of order of hundreds or thousands, but are prohibitively slow for datasets one hundred times bigger, the datasets that typically occur in the applications mentioned in Section 3.

In such a case, instead of pursuing the mathematical vision of ideal algorithm, it may be more useful to think carefully what in fact is wanted. We believe that HPŠ, unlike "certain mathematicians without personal knowledge of the subject" (to borrow an idiom from R. A. Fisher's repertory), know well that exact computation is an illusion. One can start from linear programming: for handling large datasets, *iterative* interior point methods are used—compare [9]. Even when sticking to the simplex method, hoping that the reality will exhibit its average probabilistic rather than its exponential worst-case complexity, we have to keep in mind that many of the standard, "exact" linear algebra operations are iterative, that is, approximate. And even if agreeing to use only "pure" Gaussian elimination, we still stumble upon floating-point computer arithmetic; exact algorithms in computational geometry require the use of *exact arithmetic*, an algorithmic toolbox typically not under the command of statisticians (and often neither of computer scientists).

If we want to draw a curve, we draw one hundred or one thousand small segments; if we want to solve an optimization, we employ an iterative method. For computing depth contours in KM, we do not propose to *sample* the unit sphere $\mathbb{S}^{k-1}$ (that is, we do not use any random mechanism, as the language of HPŠ might imply), but select $K$ (typically several thousand) equispaced (more sophisticated strategies are possible) directions in $\mathbb{S}^1$ (with no ambition for general $k$). Computing univariate quantiles is simple and quick—the standard R [13] function quantile uses the straightforward $O(n \log n)$ algorithm, the more sophisticated kuantile of [9] an algorithm with the best possible, $O(n)$ complexity; the ambiguities arising with possible nonunique results are solved in a straightforward fashion. (The ambiguities in linear programming solutions can be dealt with as well, but require an implementation that is returning all extremal points of the solu-tion set, not a usual feature; also, unlike with quantiles of projections, where we know that we have to use the inf version, it is not clear which version of linear-programming solution yields the right thing. Centroid? Steiner's point?) The most expensive operation here, $O(n)$, is the vector multiplication obtaining the projection itself; the contour constructing algorithm has complexity linear in $K$, as the directions come already sorted. The summary complexity is $O(nK)$, and the program runs pretty fast, easily accommodating data in the size of millions. The resulting contour is indeed approximate, but KM prove theoretically (Theorems 8 and 9) that the quality of approximation increases with $K$; it may decrease with $k$, even "extremely" as HPŠ claim (despite showing any evidence for that)—but as already mentioned, general $k$ was not our ambition.

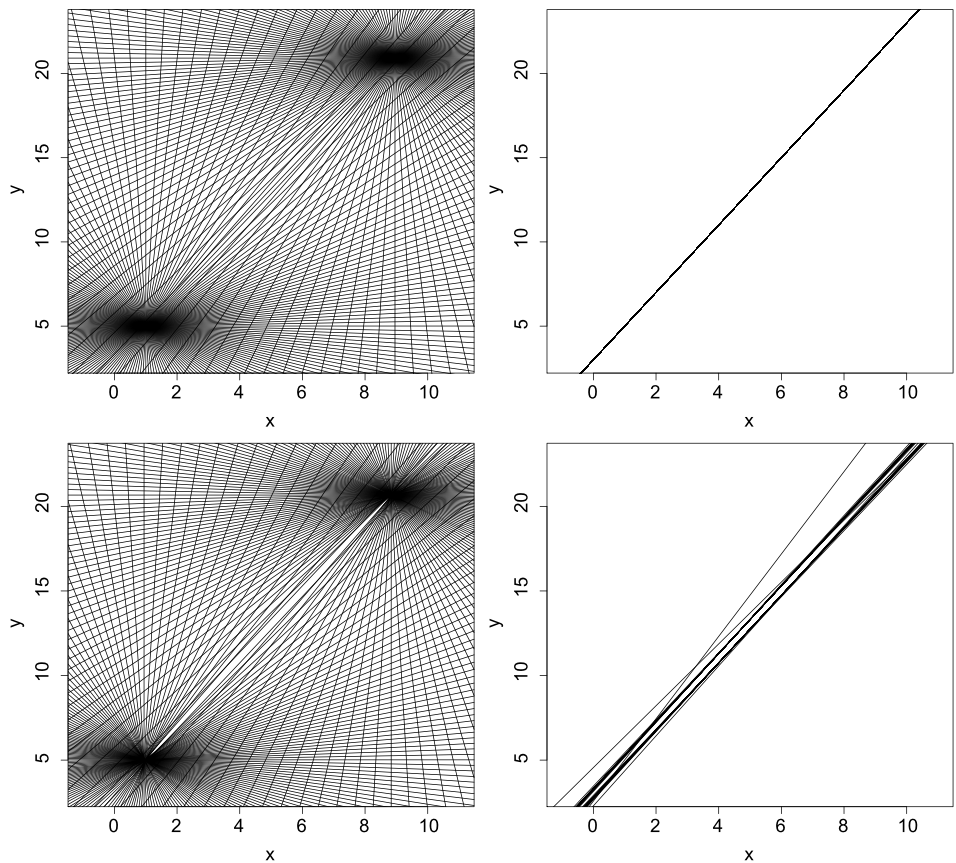

Fig. 2. *Left panels show the quantiles of projections and right the directional regression quantiles, each for* 201 *equispaced directions. The lower row shows the results for randomly perturbed data.*

Let us conclude this section with an example showing pitfalls of certain templates of thinking. Consider the dataset with points that all lie on a line (this is perhaps somehow contrived, but is the situation with perfectly linearly dependent data that impossible?). For simplicity, we will take them uniformly spaced on some segment. The depth contour for $p = 0.1$, say, is then exactly the part of the segment that contains 80% of the datapoints, with 10% on each side chopped off. If the orientation of the segment is random, the method proposed by KM is indeed unlikely to pick up the exact facet; but the result, with only 201 directions used, is pretty close (as perhaps can be seen in the upper left panel of Figure 2). On the other hand, directional quantile regressions all coincide with the line through the segment; to chop off the tails, and thus obtain the correct solution, the segment needs to be oriented vertically; all potential quantile regression then intersect in the right point, but the standard quantile regression stumbles upon the fact that the design matrix is singular, making typical quantile regression programs exit with an error message; would be interesting to watch how the parametric linear programming code of HPŠ handles that.

Nevertheless, not having the latter at hand, we *palliate* (or cure?) the problem by adding some random noise, as shown in the lower panels of Figure 2. Then both approaches return the same answer (as can be perhaps seen on the lower right panel of Figure 2).

**5. Conclusion.** We are grateful to Editors of *The Annals of Statistics* and especially to Xuming He for inviting us to write this note, and for valuable editorial suggestions; we also greatly benefited from discussions with Roger Koenker. The second author would like to return warmly the thanks addressed to his person in the Acknowledgement of the discussed paper: whatever he did, he did with pleasure.

Department of Mathematical
and Statistical Sciences
University of Alberta
CAB 632, Edmonton, AB, T6G 2G1
Canada
E-mail: lkong@math.ualberta.ca
mizera@stat.ualberta.ca